\documentclass{amsart}

\usepackage{amsmath}
\usepackage{amsfonts}
\usepackage{amssymb,enumerate}
\usepackage{amsthm}
\usepackage[all]{xy}
\usepackage{hyperref}

\newcommand{\Min}{\operatorname{Min}}

\newtheorem{thm}{Theorem}[section]
\newtheorem{cor}[thm]{Corollary}
\newtheorem{lem}[thm]{Lemma}
\newtheorem{prop}[thm]{Proposition}

\newtheorem{rem}[thm]{Remark}

\begin{document}

\bibliographystyle{amsplain}

\date{}
\author{Parviz Sahandi}

\address{Department of Mathematics, University of Tabriz, Tabriz, Iran and School of Mathematics, Institute for
Research in Fundamental Sciences (IPM), Tehran, Iran.}

\email{sahandi@tabrizu.ac.ir, sahandi@ipm.ir}



\title{Minimal prime ideals and semistar operations}

\begin{abstract}

Let $R$ be a commutative integral domain and let $\star$ be a
semistar operation of finite type on $R$, and $I$ be a
quasi-$\star$-ideal of $R$. We show that, if every minimal prime
ideal of $I$ is the radical of a $\star$-finite ideal, then the set
$\Min(I)$ of minimal prime ideals over $I$ is finite.

\end{abstract}

\maketitle

\section{Introduction}

In \cite[Theorem 88]{K}, Kaplansky proved that: Let $R$ be a
commutative ring satisfying the \emph{ascending chain condition}
(a.c.c. for short) on radical ideals, and let $I$ be an ideal of
$R$. Then there are only a finite number of prime ideals minimal
over $I$.

This result was generalized in \cite[Theorem 1.6]{GH} by showing
that (see also \cite{A}): Let $R$ be a commutative ring with
identity, and let $I\neq R$ be an ideal of R . If every prime ideal
minimal over $I$ is the radical of a finitely generated ideal, then
there are only finitely many prime ideals minimal over $I$.

In 1994, Okabe and Matsuda \cite{OM} introduced the concept of
\emph{semistar operation} to extend the notion of classical
\emph{star operations} as described in \cite[Section 32]{G}. Star
operations have been proven to be an essential tool in
\emph{multiplicative ideal theory}, allowing one to study different
classes of integral domains. Semistar operations, thanks to a higher
flexibility than star operations, permit a finer study and new
classifications of special classes of integral domains.

Throughout this note let $R$ be a commutative integral domain, with
identity and let $K$ be its quotient field.

The purpose of this note is to prove the semistar analogue of
Kaplansky's \cite[Theorem 88]{K} and Gilmer and Heinzer
\cite[Theorem 1.6]{GH} results. More precisely we prove the
following theorem.

\vspace{.05in} \noindent{\bf Theorem.} Let $\star$ be a semistar
operation of finite type on the integral domain $R$, and $I$ be a
quasi-$\star$-ideal of $R$. If every minimal prime ideal of $I$ is
the radical of a $\star$-finite ideal, then $I$ has finitely many
minimal prime ideals. \vspace{.05in}

Now we recall some definitions and properties related to semistar
operations. Let $\overline{\mathcal{F}}(R)$ denote the set of all
nonzero $R$-submodules of $K$ and let $\mathcal{F}(R)$ be the set of
all nonzero \emph{fractional} ideals of $R$, i.e.
$E\in\mathcal{F}(R)$ if $E\in\overline{\mathcal{F}}(R)$ and there
exists a nonzero $r\in R$ with $rE\subseteq R$. Let $f(R)$ be the
set of all nonzero finitely generated fractional ideals of $R$.
Then, obviously
$f(R)\subseteq\mathcal{F}(R)\subseteq\overline{\mathcal{F}}(R)$. A
semistar operation on $R$ is a map
$\star:\overline{\mathcal{F}}(R)\rightarrow\overline{\mathcal{F}}(R)$,
$E\rightarrow E^{\star}$, such that, for all $x\in K$, $x\neq 0$,
and for all $E, F\in\overline{\mathcal{F}}(R)$, the following
properties hold:
\begin{itemize}
\item [$\star_1$] $(xE)^{\star}=xE^{\star}$;
\item [$\star_2$] $E\subseteq F$ implies that $E^{\star}\subseteq
F^{\star}$;
\item [$\star_3$] $E\subseteq E^{\star}$ and
$E^{\star\star}:=(E^{\star})^{\star}=E^{\star}$,
\end{itemize}
cf. for instance \cite{OM}. Recall that, given a semistar operation
$\star$ on $R$, for all $E, F\in\overline{\mathcal{F}}(R)$, the
following basic formulas follow easily from the axioms:
\begin{itemize}
\item [(1)] $(EF)^{\star}=(E^{\star}F)^{\star}=(EF^{\star})^{\star}=(E^{\star}F^{\star})^{\star}$;
\item [(2)] $(E+F)^{\star}=(E^{\star}+F)^{\star}=(E+F^{\star})^{\star}=(E^{\star}+F^{\star})^{\star}$;
\item [(3)] $(E\cap F)^{\star}\subseteq E^{\star}\cap F^{\star}=(E^{\star}\cap F^{\star})^{\star}$, if $E\cap
F\neq 0$.
\end{itemize}
cf. for instance \cite[Proposition 5]{OM}.

A \emph{(semi)star operation} is a semistar operation that, when
restricted to $\mathcal{F}(R)$, is a star operation (in the sense of
\cite[Section 32]{G}). It is easy to see that a semistar operation
$\star$ on $R$ is a (semi)star operation if and only if
$R^{\star}=R$.

Let $\star$ be a semistar operation on the integral domain $R$. An
ideal $I$ of $R$ is called a \emph{quasi-$\star$-ideal} of $R$ if
$I^{\star}\cap R=I$. It is easy to see that for any ideal $I$ of
$R$, the ideal $I^{\star}\cap R$ is a quasi-$\star$-ideal. An ideal
is said to be a \emph{quasi-$\star$-prime}, if it is prime and a
quasi-$\star$-ideal. Let $\star$ be a semistar operation, put
$E^{\star_f}=\bigcup F^{\star}$ where the union taken over all
finitely generated $F\subseteq E$, for every
$E\in\overline{\mathcal{F}}(R)$. It is easy to see that $\star_f$
defines a semistar operation on $R$ called \emph{the semistar
operation of finite type associated to} $\star$. Note that there is
the equality $(\star_f)_f=\star_f$. A semistar operation $\star$ is
said to be \emph{of finite type} if $\star=\star_f$; in particular
$\star_f$ is of finite type. An element
$E\in\overline{\mathcal{F}}(R)$ is said \emph{$\star$-finite}, if
$E^{\star}=F^{\star}$ for some $F\in f(R)$. Note that if $I$ is an
ideal of $R$ then $I^{\star}$ is an ideal of the overring
$R^{\star}$ of $R$. Denote by $\star$-$\Min(I)$ the set of
\emph{minimal quasi-$\star$-prime} ideals over $I$. So that when
$\star=d$ is the identity semistar operation, then
$d$-$\Min(I)=\Min(I)$. It can be seen that if $I$ is an ideal of
$R$, then $\star$-$\Min(I)\subseteq\Min(I^{\star}\cap R)$. Using
\cite[Lemma 2.3(d)]{EF} we have each minimal prime over a
quasi-$\star_f$-ideal is a quasi-$\star_f$-ideal. Therefore If $I$
is a quasi-$\star_f$-ideal, we have $\star_f$-$\Min(I)=\Min(I)$.

The most widely studied (semi)star operations on $R$ have been the
identity $d_R$, $v_R$, and $t_R:=(v_R)_f$ operations, where
$E^{v_R}:=(E^{-1})^{-1}$, with $E^{-1}:=(R:E):=\{x\in K|xE\subseteq
R\}$. Our terminology and notation come from \cite{G}.

\section{Main result}

Before proving the main result of this paper, we need a lemma.

\begin{lem}\label{3} Suppose that $\star$ is a semistar operation of finite type on the integral domain $R$,
and that $I$ is a quasi-$\star$-ideal. Then $\sqrt I^{\star}\cap
R=\sqrt I$, that is $\sqrt I$ is also a quasi-$\star$-ideal of $R$.
\end{lem}

\begin{proof} Since $\sqrt I\subseteq\sqrt I^{\star}\cap R$, it is enough to
show that $\sqrt I^{\star}\cap R\subseteq\sqrt I$. Let $x\in\sqrt I^{\star}\cap R$. Then for every
$P\in\Min(I)$ we have $xR^{\star}\subseteq\sqrt I^{\star}\subseteq
P^{\star}$. Since $P$ is a quasi-$\star$-prime ideal of $R$ by
\cite[Lemma 2.3(d)]{EF}, we obtain that $x\in P$. Hence $x\in\sqrt
I$ as desired.
\end{proof}

\begin{rem} Suppose that $\star$ is a semistar operation of finite type on the integral domain $R$,
and that $I$ is a quasi-$\star$-ideal. Then $\sqrt{I^{\star}}\cap
R=\sqrt I$. Indeed suppose that $x\in\sqrt{I^{\star}}\cap R$, then
there is a positive integer $n$ such that $x^n\in I^{\star}$. Since
$I$ is a quasi-$\star$-ideal, and $x\in R$ we obtain that $x^n\in
I$. Hence $x\in\sqrt I$. Consequently we have $\sqrt{I^{\star}}\cap
R\subseteq\sqrt I\subseteq\sqrt{I^{\star}\cap
R}=\sqrt{I^{\star}}\cap R$, which gives us the desired equality.
\end{rem}

We next give the main result of this paper.

\begin{thm}\label{gh} Let $\star$ be a semistar operation of finite type on the integral domain $R$,
and $I$ be a quasi-$\star$-ideal of $R$. If every minimal prime
ideal of $I$ is the radical of a $\star$-finite ideal, then $I$ has
finitely many minimal prime ideals.
\end{thm}

\begin{proof} Note that $\sqrt I$ is a radical quasi-$\star$-ideal by Lemma \ref{3}. Hence since
we have $\Min(I)=\Min(\sqrt I)$, it is harmless to assume that $I$
is a radical quasi-$\star$-ideal.

Let $S=\{P_1\cdots P_n|\text{ each }P_i\text{ is a prime ideal
minimal over }I\}$. If for some $C=P_1\cdots P_n\in S$ we have
$C^{\star}\subseteq I^{\star}$, hence $C\subseteq C^{\star}\cap
R\subseteq I^{\star}\cap R=I$. Then any prime ideal minimal over $I$
contains some $P_i$, so $\{P_1,\cdots,P_n\}$ is the set of minimal
prime ideals of $I$. Hence we may assume that $C^{\star}\nsubseteq
I^{\star}$ for each $C\in S$. Consider the set $\mathcal{A}$
consisting of all radical quasi-$\star$-ideals $J$ of $R$ containing
$I$ such that $C^{\star}\nsubseteq J^{\star}$ for each $C\in S$.
Since $I\in \mathcal{A}$ we have $\mathcal{A}\neq\emptyset$. The set
$\mathcal{A}$ is partially ordered under inclusion $\subseteq$, and
we show that it is inductive under this ordering. Let
$\{J_{\alpha}\}_{\alpha\in\Gamma}$ be a chain in $\mathcal{A}$. Put
$J=\cup J_{\alpha}$. So that $J$ is a radical quasi-$\star$-ideal of
$R$ containing $I$ such that $J^{\star}=(\cup
J_{\alpha})^{\star}=\cup J_{\alpha}^{\star}$, in which the second
equality holds, since $\star$ is a semistar operation of finite
type. Now assume that $C^{\star}\subseteq J^{\star}$ for some $C\in
S$. Suppose that $C=P_1\cdots P_n$, and that $P_i=\sqrt{L_i}$ for
some $\star$-finite ideal $L_i$ of $R$, for $i=1,\cdots,n$. Let
$F_i$s be finitely generated ideals of $R$, such that
$L_i^{\star}=F_i^{\star}$, for $i=1,\cdots,n$. Now consider
$(F_1\cdots F_n)^{\star}=(L_1\cdots L_n)^{\star}\subseteq (P_1\cdots
P_n)^{\star}\subseteq\cup J^{\star}_{\alpha}$, which implies that
$F_1\cdots F_n\subseteq\cup J_{\alpha}$. Therefore there exists an
index $\alpha\in\Gamma$ such that $F_1\cdots F_n\subseteq
J_{\alpha}$ and hence $(L_1\cdots L_n)^{\star}=(F_1\cdots
F_n)^{\star}\subseteq J^{\star}_{\alpha}$. Thus $L_1\cdots
L_n\subseteq J_{\alpha}$, and so we have $P_1\cdots
P_n\subseteq\sqrt{L_1\cdots
L_n}\subseteq\sqrt{J_{\alpha}}=J_{\alpha}$. Consequently $(P_1\cdots
P_n)^{\star}\subseteq J^{\star}_{\alpha}$ which is impossible. Now
Zorn's Lemma gives us a maximal element $Q$ of $\mathcal{A}$. One
can actually assume that $Q\neq R$. We show that $Q$ is a prime
ideal of $R$. To this end let $a$, $b$ be two elements of $R$ such
that $ab\in Q$ and assume that $a$, $b\notin Q$. Since
$Q\varsubsetneq(Q+aR)\subseteq\sqrt{(Q+aR)^{\star}\cap R}$, and
$\sqrt{(Q+aR)^{\star}\cap R}$ is a radical quasi-$\star$-ideal (by
Lemma \ref{3}) containing $I$, there exists an element $C_1\in S$
such that $C_1^{\star}\subseteq\sqrt{(J+aR)^{\star}\cap R}^{\star}$.
By the same reason there exists again an element $C_2\in S$ such
that $C_2^{\star}\subseteq\sqrt{(Q+bR)^{\star}\cap R}^{\star}$.
Therefore we have
\begin{align*}
(C_1C_2)^{\star}=(C_1^{\star}C_2^{\star})^{\star} \subseteq & (\sqrt{(Q+aR)^{\star}\cap R}^{\star}\sqrt{(Q+bR)^{\star}\cap R}^{\star})^{\star} \\[1ex]
        = & (\sqrt{(Q+aR)^{\star}\cap R}\sqrt{(Q+bR)^{\star}\cap R})^{\star} \\[1ex]
        = & (\sqrt{((Q+aR)^{\star}\cap R)((Q+bR)^{\star}\cap R)})^{\star} \\[1ex]
        \subseteq & (\sqrt{((Q+aR)^{\star}(Q+bR)^{\star})\cap R})^{\star} \\[1ex]
        = & (\sqrt{((Q^{\star})^2+aQ^{\star}+bQ^{\star}+abR^{\star})\cap R})^{\star} \\[1ex]
        \subseteq & (\sqrt{Q^{\star}\cap R})^{\star}=(\sqrt Q)^{\star}=Q^{\star}, \\[1ex]
\end{align*}
which is a contradiction. Therefore $Q$ is a prime ideal of $R$. But
since $I\subseteq Q$, it contains a prime ideal $P$ minimal over $I$
by \cite[Theorem 10]{K}. Thus $P\in S$ and $P^{\star}\subseteq
Q^{\star}$, a contradiction.
\end{proof}

Defining different semistar operations, we can derive different
corollaries.

\begin{cor} (Gilmer and Heinzer \cite[Theorem 1.6]{GH} and Anderson
\cite{A}) Let $R$ be an integral domain, and let $I$ be an ideal of
$R$. If each minimal prime of the ideal $I$ is the radical of a
finitely generated ideal, then $I$ has only finitely many minimal
primes.
\end{cor}

The following result proved recently by El Baghdadi and Gabelli
\cite[Proposition 1.4]{EG} over P$v$MDs. They used the lattice
isomorphism between the lattice of t-ideals of R and the lattice of
ideals of the $t$-Nagata ring of $R$ over P$v$MDs.

\begin{cor} Let $I$ be a $t$-ideal of the integral domain $R$. If each minimal
prime ideal of $I$ is the radical of a $t$-finite ideal, then
$\Min(I)$ is finite.
\end{cor}

\begin{cor} Suppose that $R$ has a Noetherian overring $S$.
If $I$ is an ideal of $R$ such that $IS\cap R=I$, then $\Min(I)$ is
finite.
\end{cor}

\begin{proof} Define a semistar operation $\star$, by $E^{\star}=ES$, for each
$E\in\overline{\mathcal{F}}(R)$. Thus $I=IS\cap R=I^{\star}\cap R$
is a quasi-$\star$-ideal of $R$. Note that $\star$ is a semistar
operation of finite type. Let $P\in\Min(I)$. Since $S$ is a
Noetherian ring, we have
$P^{\star}=PS=(x_1,\cdots,x_n)S=(x_1,\cdots,x_n)^{\star}$, for some
elements $x_1,\cdots,x_n$ of $P$. This means that $P$ is a
$\star$-finite ideal. Now the result is clear from Theorem \ref{gh}.
\end{proof}

Recall that a semistar operations $\star$ on the integral domain $R$
is called \emph{stable}, if $(E\cap F)^{\star}=E^{\star}\cap
F^{\star}$ for every $E$, $F$ in $\overline{\mathcal{F}}(R)$. Again
recall from the introduction that $\star$-$\Min(I)$ is the set of
quasi-$\star$-prime ideals minimal over $I$.

\begin{lem} Suppose that $\star$ is a semistar operation stable and of finite type on the
integral domain $R$, and that $I$ is a nonzero ideal of $R$, such
that $I^{\star}\cap R\neq R$. Then
$\star$-$\Min(I)=\Min(I^{\star}\cap R)$. In particular
$$\sqrt{I^{\star}}\cap R=\sqrt{I^{\star}\cap
R}=\cap_{P\in\star\text{-}\Min(I)}P.$$
\end{lem}

\begin{proof} One sees easily that $\star$-$\Min(I)\subseteq\Min(I^{\star}\cap
R)$. For the reverse inclusion, let $P\in\Min(I^{\star}\cap R)$. So
that $I\subseteq I^{\star}\cap R\subseteq P$. Choose by
\cite[Theorem 10]{K} a prime ideal $Q$ minimal over $I$ contained in
$P$. Note that $Q$ is a quasi-$\star$-ideal, since it is contained
in $P$ (\cite[Corollary 3.9, Lemma 4.1, and Remark 4.5]{FH}). Then
$I^{\star}\subseteq Q^{\star}\subseteq P^{\star}$ and so $I\subseteq
I^{\star}\cap R\subseteq Q^{\star}\cap R=Q\subseteq P^{\star}\cap
R=P$. Thus $Q=P$, since $P$ is minimal over $I^{\star}\cap R$.
\end{proof}

\begin{cor} Let $\star$ be a semistar operation stable and of finite type on the integral
domain $R$, and $I$ be a nonzero ideal of $R$, such that
$I^{\star}\cap R\neq R$. If every quasi-$\star$-prime ideal minimal
over $I$ is the radical of a $\star$-finite ideal, then
$\star$-$\Min(I)$ is finite.
\end{cor}

\begin{proof} By the above lemma we have $\star$-$\Min(I)=\Min(I^{\star}\cap R)$.
Thus every prime ideal minimal over $I^{\star}\cap R$ is the radical
of a $\star$-finite ideal. Noting that $I^{\star}\cap R$ is a
quasi-$\star$-ideal of $R$, and using Theorem \ref{gh}, we have
$\star$-$\Min(I)$ is a finite set.
\end{proof}

In \cite[Section 3]{EFP} the authors defined and studied the
\emph{semistar Noetherian domains}, that is, domains having the
ascending chain condition on quasi semistar ideals. In \cite{P}
Picozza generalize several of the classical results that hold in
Noetherian domains to the case of semistar operations stable and of
finite type, for instance, Cohen's Theorem, primary decomposition,
principal ideal Theorem, Krull intersection Theorem, etc.

\begin{cor} (\cite[Proposition 2.4(2)]{P}) Suppose that $\star$ is a stable semistar operation of finite type on the
integral domain $R$, and that $R$ is a $\star$-Noetherian domain.
Then $\star$-$\Min(I)$ is finite for every ideal $I$ of $D$.
\end{cor}

Next we give equivalent conditions that every quasi-$\star$-prime of
$R$ is the radical of a $\star$-finite ideal.

\begin{prop}\label{acc} Let $\star$ be a semistar operation of finite type on the integral domain $R$,
the following then are equivalent:
\begin{itemize}
\item [(1)] Each quasi-$\star$-prime is the radical of a
$\star$-finite ideal.
\item [(2)] Each radical quasi-$\star$-ideal is the radical of a
$\star$-finite ideal.
\item [(3)] $R$ satisfies the a.c.c. on radical
quasi-$\star$-ideals.
\end{itemize}
\end{prop}

\begin{proof} $(1)\Rightarrow(2)$ Consider the following set. $$\mathcal{A}=\{I|I=\sqrt I,
 I=I^{\star}\cap R,\text{ which is not the radical of a
$\star$-finite ideal}\}.$$ If $\mathcal{A}\neq\emptyset$, let
$\beta=\{I_{\alpha}\}$ be a chain of elements of $\mathcal{A}$. Put
$I=\cup I_{\alpha}$. Hence $I$ is a radical quasi-$\star$-ideal of
$R$ such that $I^{\star}=(\cup I_{\alpha})^{\star}=\cup
I_{\alpha}^{\star}$. Suppose that $I=\sqrt L$ for some
$\star$-finite ideal $L$. Let $L^{\star}=F^{\star}$ for some $F\in
f(R)$. So that $F^{\star}=L^{\star}\subseteq I^{\star}=\cup
I_{\alpha}^{\star}$ which implies that $F\subseteq\cup I_{\alpha}$.
Therefore there is an index $\alpha$ such that $F\subseteq
I_{\alpha}$. Consequently $L^{\star}=F^{\star}\subseteq
I_{\alpha}^{\star}$ and hence $L\subseteq I_{\alpha}$. So we obtain
that $I_{\alpha}=\sqrt L$, which is impossible. Hence
$I\in\mathcal{A}$. Thus by Zorn's Lemma $\mathcal{A}$ has a maximal
element $P$. Let $a$, $b$ be two elements of $R$ such that $ab\in P$
and suppose that $a$, $b\notin P$. Since
$P\varsubsetneq(P+aR)\subseteq\sqrt{(P+aR)^{\star}\cap R}$, and
$\sqrt{(P+aR)^{\star}\cap R}$ is a radical quasi-$\star$-ideal (by
Lemma \ref{3}), we have $\sqrt{(P+aR)^{\star}\cap R}=\sqrt L$, for
some $\star$-finite ideal. By the same reason
$\sqrt{(P+bR)^{\star}\cap R}=\sqrt N$, where $N$ is a $\star$-finite
ideal. The same proof as Theorem \ref{gh} shows that $P=\sqrt{LN}$,
which is impossible, since $LN$ is a $\star$-finite ideal. Hence $P$
is a quasi-$\star$-prime, a contradiction. Hence
$\mathcal{A}=\emptyset$.

$(2)\Rightarrow(3)$ Suppose that $(I_n)_{n\in \mathbb{N}}$  be an
ascending chain of radical quasi-$\star$-ideals, and set
$I=\bigcup_{n\in \mathbb{N}} I_n$. Then $I$ is a radical
quasi-$\star$-ideal. Hence $I=\sqrt L$ for some $\star$-finite ideal
$L$. So, there is an integer $n_0$ such that $I_{n_0}=\sqrt L=I$.
Hence $(I_n)_{n\in \mathbb{N}}$ is stationary.

$(3)\Rightarrow(1)$ Suppose that $(1)$ is false. Then we can
construct a chain $(I_n)_{n\in \mathbb{N}}$ of radical
quasi-$\star$-ideals strictly ascending. Indeed, let $P$ be a
quasi-$\star$-prime ideal which is not the radical of a
$\star$-finite ideal. Set $I_1=\sqrt{(x)}$, where $0\neq x\in P$.
Given $I_n=\sqrt{(x_1,\cdots,x_n)^{\star}\cap R}$, where
$x_1,\cdots,x_n\in P$, then
$I_{n+1}=\sqrt{(x_1,\cdots,x_n,x_{n+1})^{\star}\cap R}$, where
$x_{n+1}\in P\backslash I_n$.
\end{proof}

\begin{cor}\label{K} Suppose that $\star$ is a
semistar operation of finite type on the integral domain $R$. If $R$
satisfies the a.c.c. on radical quasi-$\star$-ideals, then
$\star$-$\Min(I)$ is finite for every ideal $I$ of $R$.
\end{cor}

\begin{proof} Note $\star$-$\Min(I)\subseteq\Min(I^{\star}\cap R)$ and use Theorem
\ref{gh}.
\end{proof}

\begin{rem} (1) One can prove Theorem \ref{gh} for arbitrary rings with zero divisors.
Let $R$ be a commutative ring, with total quotient ring $T(R)$. Let
$\mathcal{F}(R)$ denote the set of all $R$-submodules of $T(R)$.
Suppose an operation $*:\mathcal{F}(R)\rightarrow\mathcal{F}(R)$,
$E\rightarrow E^*$, satisfies, for all $E, F\in\mathcal{F}(R)$, and
for all $x\in T(R)$, the following:
\begin{itemize}
\item [$*_1$] $xE^*\subseteq (xE)^*$ and if $x$ is regular, then $xE^*=(xE)^*$;
\item [$*_2$] $E\subseteq F$ implies that $E^*\subseteq F^*$;
\item [$*_3$] $E\subseteq E^*$ and
$E^{**}:=(E^*)^*=E^*$.
\end{itemize}
Then from these axioms the following directly follow:
\begin{itemize}
\item [(i)] $(EF)^*=(E^*F)^*=(EF^*)^*=(E^*F^*)^*$;
\item [(ii)] $(E+F)^*=(E^*+F)^*=(E+F^*)^*=(E^*+F^*)^*$;
\item [(iii)] $(E\cap F)^*\subseteq E^*\cap F^*=(E^*\cap
F^*)^*$.
\end{itemize}

It is clear that any semistar operation satisfies these axioms.

It is routine to see that \cite[Lemma 3.3]{HP} the $v$-operation
satisfies, these axioms, where $E^v=(E^{-1})^{-1}$, in which
$E^{-1}=(R:E)=\{x\in T(R)|xE\subseteq R\}$, for $E\in
\mathcal{F}(R)$.

By this operation Theorem \ref{gh}, is true for rings with zero
divisors.

%Let $I$ be an ideal of the ring $R$, and $\star$ be a semistar
%operation on it. An ideal $I$ of $R$ is called a
%\emph{quasi-$\star$-ideal} of $R$ if $I^{\star}\cap R=I$. An ideal
%is said to be a \emph{quasi-$\star$-prime}, if it is prime and
%quasi-$\star$-ideal. Recall that $\star$ is of finite type if for
%every ideal I of R; $I^{\star}=\sum\{J^{\star}|\text{ J is a
%finitely generated ideal of R contained in }I\}$. An element
%$E\in\mathcal{F}(R)$ is said \emph{$\star$-finite}, if
%$E^{\star}=F^{\star}$ for some $F\in f(R)$. As in the domain case
%\cite[Lemma 2.3(d)]{EF}, any prime ideal minimal over a
%quasi-$\star$-ideal is a quasi-$\star$-ideal provided that $\star$
%is of finite type.

(2) It is interesting to note that if we take $R$ to be the ring of
all sequences from $\mathbf{Z}/2\mathbf{Z}$ that are eventually
constant, with pointwise addition and multiplication, then $R$ is a
zero-dimensional Boolean ring with minimal prime ideals
$P_i=\{\{a_n\}\in R|a_i=0\}$ and $P_{\infty}=\{\{a_n\}\in
R|a_n=0\text{ for large }n\}$ and each $P_i$ is principal but
$P_{\infty}$ is not finitely generated. Thus while $R$ has
infinitely many minimal prime ideals, only one is not the radical of
a finitely generated ideal.
\end{rem}

In the rest of the paper we will define a class of rings, that
satisfy the a.c.c. on radical quasi-$\star$-ideals.

Let $R$ be a commutative ring. An ideal $I$ of $R$ is called an
\emph{ideal of strong finite type} (\emph{SFT-ideal} for short) if
there exist a finitely generated ideal $J\subseteq I$ and a positive
integer $k$ such that $a^k\in J$ for each $a\in I$. The ring $R$ is
called an \emph{SFT-ring} if each ideal of $R$ is an SFT-ideal.
These concepts were introduced by J. T. Arnold in \cite{Ar}. The
condition that $R$ is an SFT-ring plays a key role in computing the
Krull dimension of the power series ring $R[[X]]$ over $R$. In
\cite{K1}, Kang and Park defined and studied the $\star=t$ analogue
of SFT-rings. Now we define the more general semistar-SFT-rings.

Let $R$ be a domain and $\star$ a semistar operation on it. We
define a nonzero ideal $I$ of $R$ to be a \emph{$\star$-SFT-ideal}
if there exist a finitely generated ideal $J\subseteq I$ and a
positive integer $k$ such that $a^k\in J^{\star}$ for each $a\in
I^{\star_f}$. The ring $R$ is said to be a \emph{$\star$-SFT-ring}
if each nonzero ideal of $R$ is a $\star$-SFT-ideal. Obvious
examples of a $\star$-SFT-ring is $\star$-Noetherian domains.

\begin{prop} Suppose that $\star$ is a semistar operation of finite type on the
integral domain $R$. If $R$ is a $\star$-SFT-ring, then $R$
satisfies the a.c.c. on radical quasi-$\star$-ideals.
\end{prop}

\begin{proof} Let $P$ be a quasi-$\star$-prime ideal. Since $P$ is a $\star$-SFT-ideal
there is a finitely generated subideal $J\subseteq P$ such that
$\sqrt{P^{\star}}=\sqrt{J^{\star}}$. Now consider
$$P=\sqrt P=\sqrt{P^{\star}\cap R}=\sqrt{P^{\star}}\cap R=\sqrt{J^{\star}}\cap
R=\sqrt{J^{\star}\cap R}.$$ Since $J^{\star}\cap R$ is an
$\star$-finite ideal, the result follows by Proposition \ref{acc}.
\end{proof}

\begin{cor} Each quasi-$\star$-ideal of a $\star$-SFT-ring $R$, has only finitely many minimal prime ideals.
\end{cor}

We close the paper with the following characterization of
$\star$-SFT-rings.

\begin{prop} Suppose that $\star$ is a semistar operation of finite type on the
integral domain $R$. Then $R$ is a $\star$-SFT-ring if and only if
each quasi-$\star$-prime ideal of $R$ is a $\star$-SFT-ideal.
\end{prop}

\begin{proof} Suppose that $R$ is not a $\star$-SFT-ring. Therefore the set
$$\mathcal{A}=\{I|I=I^{\star}\cap R,\text{ and is not a
}\star\text{-SFT-ideal}\},$$ is not an empty set. The set
$\mathcal{A}$ is partially ordered under inclusion, and is inductive
under this ordering. By Zorn's lemma, $\mathcal{A}$ contains a
maximal element $P$. Assume that $a_1, a_2$ be two elements of $R$
such that $a_1a_2\in P$ and $a_1, a_2\notin P$. Since
$P\varsubsetneq (P+a_iR)^{\star}\cap R$, $(P+a_iR)^{\star}\cap R$ is
a $\star$-SFT-ideal of $R$. Consequently there exist a finitely
generated ideal $L_i\subseteq (P+a_iR)^{\star}\cap R$, and a
positive integer $k_i$ such that $c^{k_i}\in (L_i)^{\star}$ for each
$c\in (P+a_iR)^{\star}$. Let $L=L_1L_2$ and $k=k_1+k_2$. Then $L$ is
a finitely generated subideal of $P$ such that
$c^k=c^{k_1}c^{k_2}\in
(L_1)^{\star}(L_1)^{\star}\subseteq(L_1L_2)^{\star}$, for each $c\in
P^{\star}$. Thus $P$ is a $\star$-SFT-ideal, a contradiction.
Therefore, $P$ is a quasi-$\star$-prime ideal which is not a
$\star$-SFT-ideal.
\end{proof}

\begin{center} {\bf ACKNOWLEDGMENT}

\end{center}
 I would like to thank Professor Marco Fontana for his useful
comments on this paper.

%\begin{center} {\bf ACKNOWLEDGMENT}

%\end{center}

\end{document}